# Time-Spatial Serials Differences' Probability Distribution of Natural Dynamical Systems

Wei Ping Cheng[1]*, Zhi Hong Zhang[2], Pu Wang[3]

[1]*College of Civil Engineering and Architecture, Zhejiang University, Hangzhou 310058, P. R. China.*

[2]*Department of Civil Engineering, Shanghai Normal University, Shanghai 200233, P. R. China.*

[3]*Jiangxi Xinyu Guoke Technology CO.，LTD, XinYu 338000, P. R. China.*

Corresponding author: Wei Ping Cheng

**Email:** chengweiping@zju.edu.cn

**Abstract**

The normal distribution is used as a unified probability distribution, however, our researcher found that it is not good agreed with the real-life dynamical system's data. We collected and analyzed representative naturally occurring data series (e.g., the earth environment, sunspots, brain waves, electrocardiograms, some cases are classic chaos systems and social activities). It is found that the probability density functions (PDFs) of first or higher order differences for these datasets are consistently fat-tailed bell-shaped curves, and their associated cumulative distribution functions (CDFs) are consistently *S-shaped* when compared to the near-straight line of the normal distribution CDF. It is proved that this profile is not because of numerical or measure error, and the *t*-distribution is a good approximation. This kind of PDF/CDF is a universal phenomenon for independent time and space series data, which will make researchers to reconsider some hypotheses about stochastic dynamical models such as Wiener process, and therefore merits investigation.

**Significance Statement**

In most condition, researchers develop their models based on the normal distribution. However, a lot of case of natural dynamical systems in this paper show that their time-spatial serials differences' probability distribution is deviated from normal distribution clearly. None of the cases are relevant. At the same time, their profiles look like a universal shape. We have test more than global 4000 hydrological stations' run-off, sunspots, earth elevation, brain waves, stock price, etc. Those cases cover planetary sciences, engineering, environmental sciences, biophysics, economic sciences, which suggest that it may be an ultimate unified description more than the normal distribution for the real-life dynamical system, not just hypothesis. Except economy, other is found firstly. Thus, a lot of stochastic analysis based on Wiener process hypothesis will need improvement, and many Monte-Carlo simulation need re-examinations too.

## Origins of questions

We observed an interesting phenomenon during our research into detecting pipe bursts in a city water supply system, where flow fluctuates over time. The fluctuations form a stationary random process (Fig.1a). We calculated the cumulative distribution function (CDF) of water pressure and the first difference of flow rate for different sets of water pressure and flow data; an example is shown in Fig.1b (*1*). What we found interesting was that all the curves had characteristically similar shapes: the middle section of the CDF curve (i.e., close to the mean value) matches a normal distribution, but the other two sections deviate significantly from the straight line, forming an asymmetric S-shaped curve. This phenomenon was observed in the flow datasets for all the Chinese cities in which we investigated pipe bursts.

We were interested by this common water flow pattern in various city water supply systems. We collected data for other dynamical systems from different natural or social systems, such as brain waves, seismic waves, sunspots, and stock markets. The CDFs of their first or higher order differences were all similar in shape to the water supply system CDFs. This observation caused us to ask the following two questions:

(1) Is this phenomenon caused by errors in instrumentation or in data processing?

(2) Is this natural or social phenomenon occasional or universal?

Question 1 is explained in Supporting Information: Bathtub effects of differential and difference of dynamical system error. The phenomenon is not caused by data processing; time step size does not affect the PDF over a limited range. Question 2 is discussed in detail in the next section. Despite differences in spatial or temporal scales and frequencies, the CDFs of first or higher order differences are all similar. We think that this statistical pattern may represent a universal natural and social phenomenon rather than an occasional occurrence.

**Probability distribution of natural and social time or space data series differences**

Definition of time series differences: three different definitions of the first order difference of a time series are given by:

$$\Delta \hat{f}(i) = \hat{f}(i) - \hat{f}(i-1) \tag{1}$$

$$\Delta \hat{f}_m(i) = \frac{\hat{f}(i) - \hat{f}(i-1)}{\sum_{i-K+1}^{i} \hat{f}(i)/K} \tag{2}$$

$$\Delta \ln \hat{f}(i) = \ln\left(\hat{f}(i)\right) - \ln\left(\hat{f}(i-1)\right) \tag{3}$$

where $\hat{f}$ is a time series, $\Delta \hat{f}(i)$ is the first order difference, $\Delta \hat{f}_m(i)$ is the difference in the ratio of the time series to the mean value of the nearest $K$ samples, and $\Delta \ln \hat{f}(i)$ is the logarithmic first order difference. These three forms all depend on short-term variation in the time series data, and each is an approximation of the first order difference. Equation (1) can be used when the sample does not greatly vary. Equation (2) is recommended if the time series (e.g., surface runoff) changes greatly or contains negative numbers because it makes the data dimensionless and thus enables the aggregation of different time series into one dataset. Equation (3) can be used when a positive variable changes to a large degree (e.g., a stock index fluctuation rate is auto-correlated). Successive second or higher order differences can be calculated from Equations (1–3). Consistency of statistical characteristics can be ensured by the choice of suitable time steps (details are given in Supporting Information). Note that any numerical errors in difference calculation from a lower-order to a higher-order difference will be amplified.

In most condition, wiener process is used as basic hypothesis in natural or social systems, which assumed Brownian motion is a normal distribution, a standard solution of the linear diffusion equation. However, our investigation does not support the hypothesis. the CDFs of 21 cases is in six categories (Table 1): (1) numerical solutions of classic chaotic systems; (2) time series data for various natural phenomena on Earth; (3) signals from artificial structures and experimental devices; (4) biological signals; (5) time series data of social behaviour; and (6) spatial differences. Most of the data samples were downloaded from public websites; their physical characteristics, sampling frequency and the measuring equipment used to obtain them are irrelevant.

**Examples of chaotic dynamical systems**: chaotic systems are neither periodic nor convergent due to their sensitivity to initial values. The Lorenz system is a set of differential equations for air flow; the Duffing equation (*2*) is a nonlinear vibration model; the Chua circuit (3) is a nonlinear electronic circuit. All of them behave chaotically. Numerical models using these equations(Fig.2) conform to our hypothesis in that at some order of difference they show a bell-shaped histogram that gives an *S-shaped* CDF. When the order of the difference increases (to an order higher than that of current variables), CDFs for these variables have similar shapes(Fig.2). The shape of the PDF is independent of the initial values and perturbations.

**Natural environmental signal in the Earth:** Phenomena such as air temperature, wind speed, tidal currents, and runoff usually vary at different temporal or spatial scales. The first case we examine is the daily temperature in Nanjing recorded at 11:00 or 11:25 from 2001-01-01 to 2018-05-25. The histogram of the temperature data contains multiple peaks (Fig.3A(a2)). In the first order difference, a single peak remains and the histogram shifts slightly to the left (cooling occurs quickly, but heating is slow). The second order difference curve becomes symmetric (Figs.3A (a4, a6)). Wind speed changes randomly at the atmospheric boundary layer. The wind speed data were recorded at Yueqing, China. The wind speed PDF is asymmetric, but the first order difference PDF is symmetric, and the corresponding CDF deviates from the central straight line at both ends (Figs.3B (b2, b3)). Tides are due to the gravitational attraction of Earth, moon and sun; they are periodic. Tide data were recorded for the East China Sea (Fig.3C). Daily surface runoff data for the River Thames in the United Kingdom were recorded from 1883-01-01 to 2017-09-20 (Fig.3; original data can be downloaded from NRFA (**4**). A large change in surface runoff data can be observed. Runoff values in flood periods can be hundreds of times those in a normal period. We used Equation (2) to investigate the small fluctuations in surface runoff data. Fig.3D shows daily runoff and daily fluctuation rate. We note that the range of daily fluctuation rate has greatly increased in the past 30–40 years. We have tested data from more than global 4000 hydrological stations, and they all have the same pattern.

**Seismic waves:** a seismic wave is a short-term time series. The three sets of seismic waves used in this study have different power spectra but they have similar first order difference CDFs (Fig. 4). Note that the seismic wave data for the first order difference of the ground velocity were measured by the accelerometer sensors (Fig. 4 D).

**Sunspots:** sunspots are comparatively dark areas on the surface of the sun where magnetic fields gather. The main cycle of sunspot intensity is about 11 years. The time series is chaotic (*5,6*). Individual sunspots or groups of sunspots may occur anywhere on the surface of the sun and remain for a period of a few days to a few months. Sunspots eventually decay and disappear. Sunspot data observed from 1977-01-01 to 2013-12-31 were downloaded from NASA (*7*). The CDF of the data is S-shaped, as shown in Fig.5.

**Artificial structures and experimental devices:** fluctuations in water pressure and flow rate in a water supply system can be due to changes in water consumption, changes in water supply controls, and measurement errors from monitoring instruments. We selected two dynamical systems that include artificial structures (a high-speed railway in China and a wind tunnel) which are unrelated to human behaviour. The data for the acceleration of ground vibrations along the high-speed railway and wind pressure data for an atmospheric boundary layer from a wind tunnel simulation have similarly shaped CDFs. The ground vibration data had a sampling frequency of 1280Hz. The 2004 wind pressure experiment was conducted at the Wind Engineering Research Center of Tokyo Polytechnic University (*8*) (Fig.6).

**Bioelectrical signals:** Electrocardiograms (ECGs) and electroencephalograms (EEGs) represent bioelectrical signals from the heart and the brain. An ECG measures change in the potential of cardiomyocytes from the endocardial–pericardial sequence. Human ECG data (*9-13*) was downloaded from PhysioNet (*14*), and an EEG of a dog (*Canis lupus familiaris*) was provided by the School of Pharmaceutical Sciences, Zhejiang University. CDFs of the ECG first-order difference deviate from the straight line (normal distribution) outside the range of ±2~3 times the variance. The CDFs have the same shape as those for other natural time series. The EEGs were downloaded from EPFL (*15*) (one sample is a patient EEG, the other is a healthy human EEG (*16*) and PhysioNet (*14*). Although EEGs are complex and noisy, CDFs of their first-order and higher order differences are similar in shape (Fig.7B(b)). A phonocardiogram (PCG) represents cardiac sounds caused by the movement of blood in the heart and blood vessels. PCG data were downloaded from PhysioNet (*17,18*). The CDF of PCG first-order difference had a similar shape to the other CDFs discussed in this paragraph.

**Social phenomena:** stock trading is a social behaviour that can easily be accurately digitized. The Dow Jones Industrial Average index (DJIA) is the first weighted average price index to have been created, in 1884. Since then it has increased by a factor of 500. Equation (3) was used to process the DJIA data. Even though the history of the Chinese stock market is shorter than the lifespan of the DJIA, it has the same profile. Each stock price on the Chinese market is a different time series. Equation (2) converts the time series of different individual stocks into dimensionless data which we then aggregated into one dataset. Figs.8B (b1, b2) are the dimensionless daily stock price data, with fluctuations, in China. Figs.8B (b3, b4) are 5-minute and 10-minute price data of 900 stocks (the parameter $K = 5$). Both the New York and Chinese stock markets, having different time scales, have stock price CDFs that deviate from the normal distribution at both ends (Figs.8A(a4), 8B(b2-b4)). The CDF for copper spot prices published by the Chilean government (*19*) has a similar shape (Fig.8).

**Spatial differences:** the data series considered so far were time series. Spatial data series contain data points that vary over multidimensional space rather than time. We collected two examples of spatial data series: one is terrestrial elevation, the other is the cosmic microwave background (CMB). Ground elevation data for 38 °N, 102 °W was downloaded from the American Geographic Society (*20*). The CDFs for both north-south and east-west directions have the same shape as those for the time series datasets (Figs.9A(a1-a4)). CMB is electromagnetic radiation that is a remnant from an early stage of the universe (*21*). The cosmic background radiation temperature map was generated using data collected by the European Space Agency (ESA) Planck satellite (*22*). Four products (SMICA, NILC, SEVEM and COMMANDER-Ruler) can been downloaded from ESA (*23*). The CDFs (Fig.9B(b3)) were generated from the SMICA product, with a resolution of 2048. The CDFs from the calculated first order differences in both directions (Fig.9B(b2)) gradually deviate from the normal distribution curve at points close to 1% or 99%, and turn sharply near 0.3% and 99.7%.

**Normalized CFD of data series differences**

Statistical curve shape of first order differences for time series or spatial data series: our analysis of the preceding cases provided three consistent characteristics of CDFs for first (or higher) order differences:

(1) The mean value of the difference is zero or close to zero; the ratio of mean value to standard deviation $\frac{\mu}{\sigma}$ can be used as an index (except for DJIA, which has been increasing over time).

(2) The PDF of first or higher order difference is basically a bell-shaped symmetric curve.

(3) The central sections of the CDF curves coincide with the normal distribution CDF, but both ends deviate significantly. The CDF curves are S-shaped and are similar to a *t*-distribution or a Cauchy distribution.

The first characteristic is easy to understand. In a finite system, over a long period, $\lim_{n \to \infty} \frac{1}{n} \sum \Delta \hat{y}(i) = 0$. The second and third characteristics have not yet been mathematically explained. Table 2 shows that curve fitting results match a *t*-distribution. All the correlation coefficients are >0.9, which indicates that a *t*-distribution is a good approximation(Fig.10).

## Discussion

We collected some typical independent time series and space series data that occur in nature, in social activity, and in artificial structures and analyzed them. Some datasets were chaotic and unstable in the conventional sense. Despite the diversity of their natural properties and their different frequency spectra, the shapes of the CDFs of first or higher order differences were similar in that the mean value was close to zero. They were also similar in that the PDF had a single peak and was a symmetric bell-shaped curve that approximated a *t*-distribution. Our statistical analysis showed that the CDF shape is universally similar.

We also discovered some very interesting phenomena. For example, the fluctuation in runoff into the River Thames has increased over the last few decades. We would like to know whether this increase is due to measurement error, human activity, or climate change. Our results show that the first order differences of a normal person ECG or a patient ECG are not white noise. These illustrative observations indicate that our analysis of first and higher order differences can establish previously unrecognized regularities in chaotic physical or social systems. There is historical precedent for this claim: for example, Newton discovered the law of momentum by examining difference in velocities; Faraday discovered the laws of electromagnetism by analyzing the rate of change of magnetic flux. Further research into the nature of first and higher order difference may uncover more physical or social regularities and lead to the recognition of underlying law. The normal distribution is widely used in science and engineering.

However, the solution of a linear diffusion equation may not be the best way to explain a nonlinear phenomenon: a fat-tailed PDF may be a more realistic representation than the normal distribution; many Monte Carlo simulations with a *t*-distribution hypothesis may be more accurate than a normal distribution.

## Acknowledgments

We wish to thank all data provider including person and organizations.

## Data availability

The authors declare that all relevant data supporting the finding of this study are available within the paper and its Supplementary Information files. Additional data are available from the corresponding author upon request. Most of the data comes from public websites. We strongly recommend that readers download it directly from the websites marked within references.

# Figures and Tables

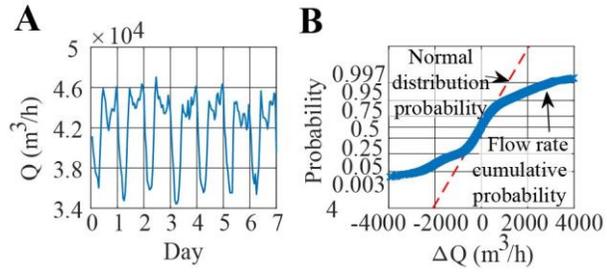

**Fig.1 Flow fluctuation in water supply system**. (**A**) flow rate over a week. (**B**) CDF of the first order difference of flow rate.

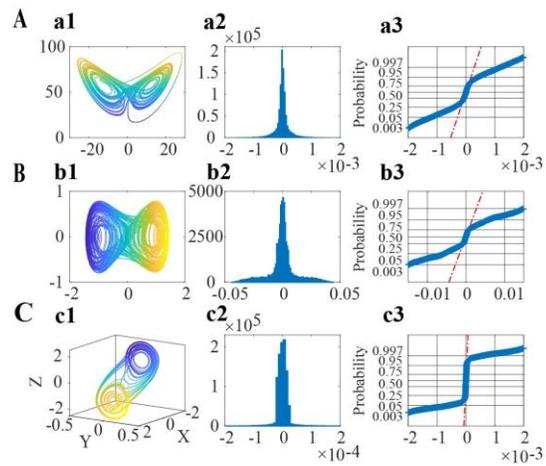

**Fig.2 Chaotic dynamical systems. (A)** Lorenz system: (a1) state space in *x–z* plane, (a2) histogram of *x* third order difference, (a3) CDF of *x* third order difference. (**B**) Duffing equation: (b1) state space in *x–z* plane, (b2) histogram of *x* second order difference, (b3) CDF of *x* second order difference. (**C**) Chua circuit: (c1) 3D graph of Chua circuit, (c2) histogram of *x* third order difference, (c3) CDF of *x* third order difference.

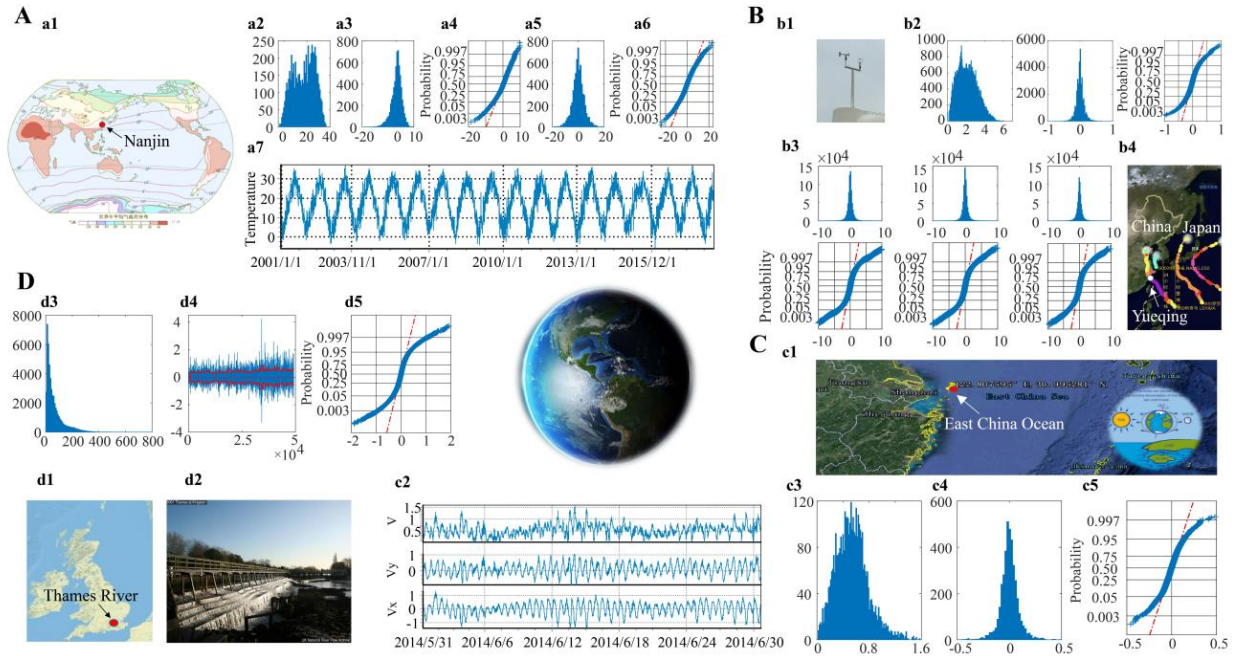

**Fig.3 Representations of data for natural geophysical Earth phenomena. (A)** characteristic statistics of noontime air temperature at Nanjing, China (unit: ℃): (a1) location of Nanjing, (a2) air temperature histogram, (a3) air temperature first order difference histogram, (a4) CDF of air temperature first order difference, (a5) air temperature second order difference histogram, (a6) CDF of air temperature second order difference. **(B)** characteristic statistics of wind speed (sampling period 0.1s, unit: m/s): (b1) ultrasonic anemometer at the top of a corner column in Yueqing Stadium, (b2) histogram of wind speed, histogram of wind speed first-order difference, and CDF of wind speed first order difference in Yueqing stadium from 11:00 to 12:00 on 2013-07-14 (wind speed was in the range 1–2 m/s), (b3) histogram and CDF of wind speed first-order difference in 3 directions, (b4) path of supertyphoon Lekima and monitoring position from 2019-09-09 14:00 to 2019-09-10 10:00 (maximum wind speed 25 m/s). **(C)** tidal current velocity in the East China Sea (30.995281°N, 122.807595°E) from 2014-05-31 to 2014-06-30 (sampling interval 10 minutes, unit: m/s): (c1) tidal current monitoring position, (c2) tidal current speed, (c3) histogram of current speed, (c4) histogram of current speed first order difference, (c5) CDF of current speed first order difference. **(D)** characteristic statistics of daily runoff into the River Thames over 135 years (unit: m$^3$/d): (d1) monitoring position, (d2) photo of monitoring position, (d3) histogram of daily runoff, (d4) daily runoff first order difference (the red line is the standard deviation for daily flow over two years; it has increased by 30%–50% in the past 30–40 years), (d5) CDF of daily runoff first order difference.

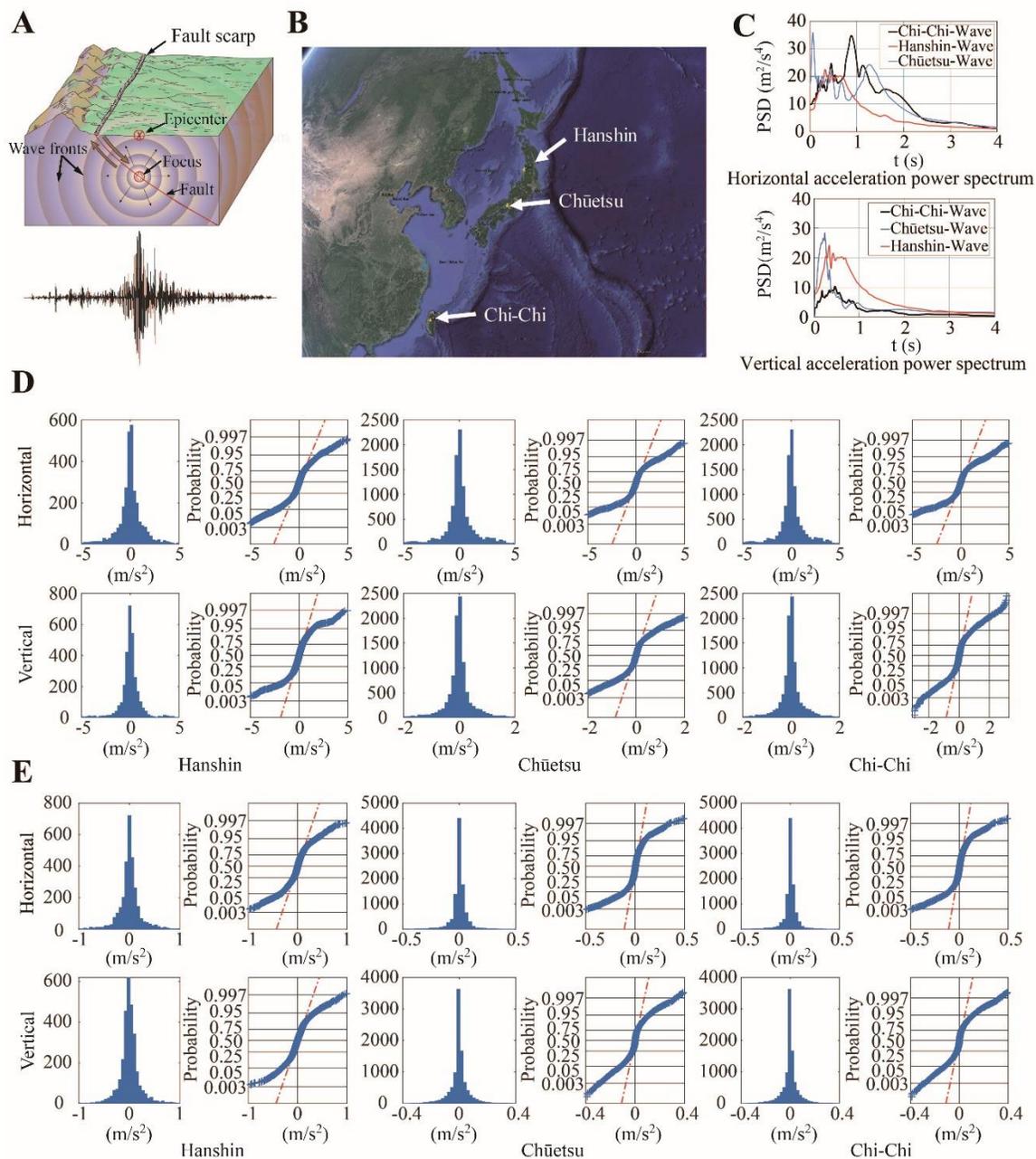

**Fig.4. Seismic waves (principal earthquake wave).** (**A**) seismic wave propagation and monitoring. (**B**) locations of seismic wave sources: Hanshin-Awaji and Chūetsu earthquakes in Japan, and Chi-Chi earthquake in Taiwan (the sampling period of Hanshin and Chūetsu is 0.01 s and of Chi-Chi is 0.005 s). (**C**) spectra of earthquake wave (there is no relationship between them; even in the same earthquake. the horizontal and vertical power spectra are quite different). (**D**) histogram of earthquake acceleration and CDF (the seismic accelerometer had calculated the first order difference of ground speed). (**E**) histogram and CDF of earthquake acceleration first order difference.

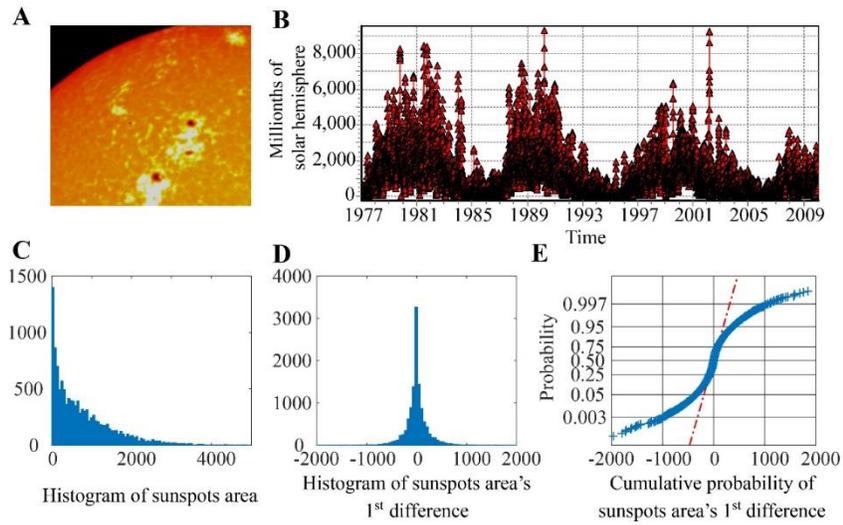

**Fig.5. Sunspot statistical distributions. (A)**, sunspots; **(B)** daily area of sunspots from 1977 to 2013. **(C)** sunspot area histogram. **(D)** histogram of sunspot area first order difference. **(E)** CDF of sunspot area first order difference.

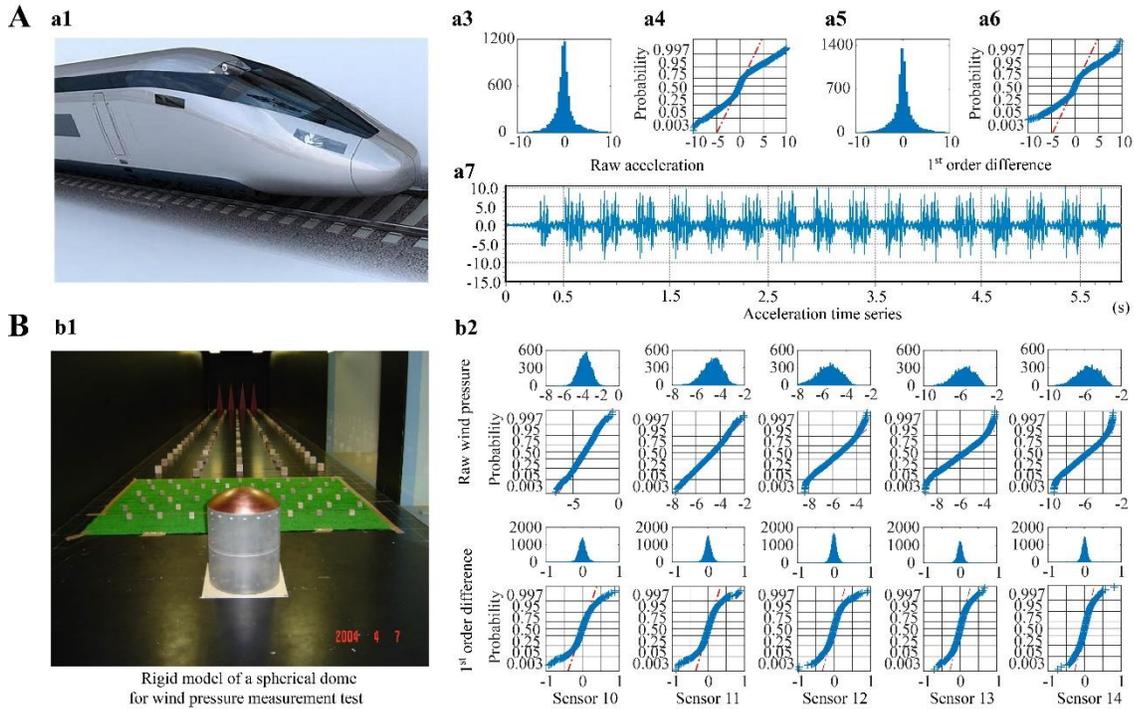

**Fig.6. Artificial structures and experiments. (A)** acceleration signal of ground vibration along a high-speed railway(unit: mv): (a1) high speed railway data were obtained for the up line of the Shanghai–Nanjing high speed railway in 2016, at a speed of 200–300 km/h; (a2) serial acceleration data serial with sampling frequency 1280 Hz (the vibration sensor uses a DC differential wiring method to filter out low frequency components below 0.7 Hz); (a3) histogram of acceleration; (a4) CDF of acceleration data; (a5) histogram of acceleration first order difference; (a6) CDF of data first-order difference. **(B)** wind pressure measurement in a wind tunnel (unit: Pa): (b1) wind tunnel test in the wind tunnel (which has an operational section 15 m long, 2.2 m wide and 1.8 m high) of the Wind Engineering Research Center of Tokyo Polytechnic University: a thickness $Z_G$ = 100 cm of two boundary layers with terrain types having power laws 0.195 and 0.346 was generated using a set of spikes and a number of roughness blocks on the wind tunnel floor using a sampling interval of 15 s; (b2) histograms of wind pressure and CDFs of wind pressure first order difference.
18

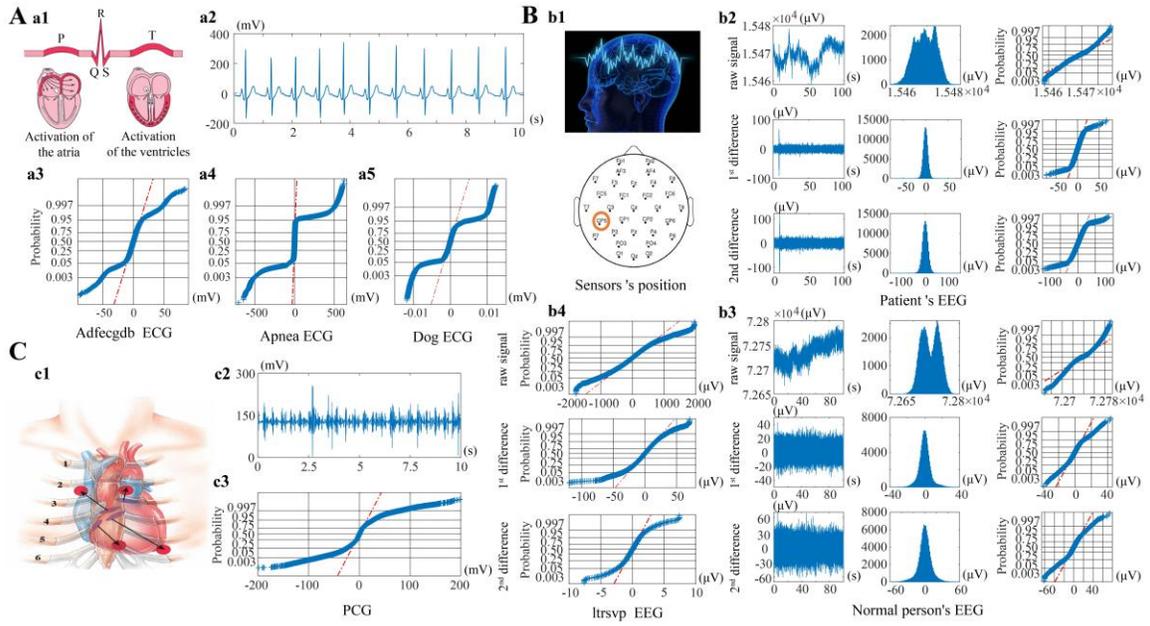

**Fig.7 Bioelectrical signals (ECG, EEG and PCG, some images from internet).** (**A**) electrocardiogram, (a2) electrocardiogram from Apnea–ECG Database, (a3, a4) CDF of human ECG first order difference, (a5) CDF of dog ECG first order difference. (**B**) brain wave diagram, (b2, b3) raw time series data, CDFs and histograms of first and second order differences, (b4) ltrsvp EEG at PO8, including raw signals, first and second order difference data, and CDFs. (**C**) cardiac sound diagram, (c2) PCG time series source, (c3) CDF of PCG first order difference.

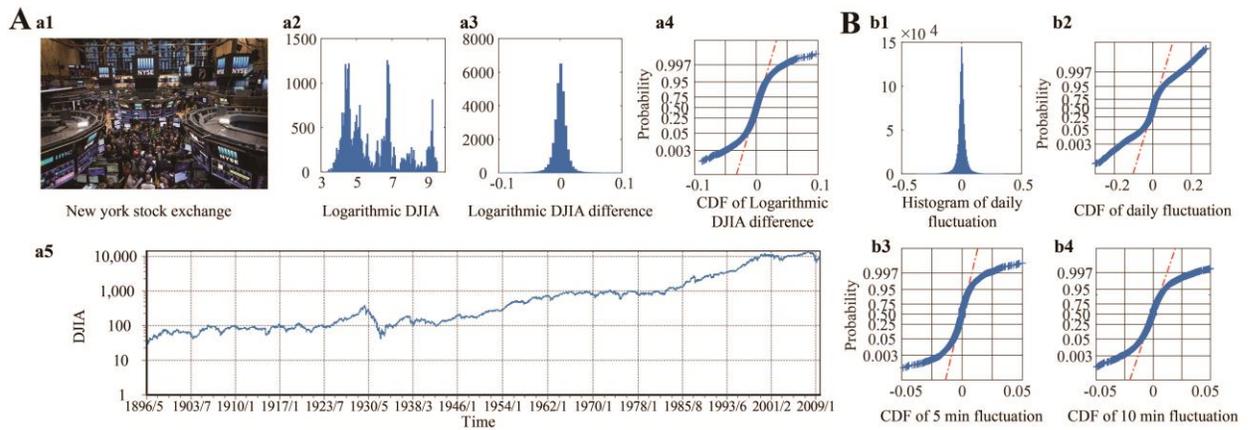

**Fig.8 Stock market fluctuations. (A)** Dow Jones Industrial Average index: (a1) New York Stock Exchange (image from internet), (a2) logarithmic DJIA histogram, (a3) DJIA histogram and CDF of logarithmic first order difference. (**B**) China stock market stock prices: (b1) histogram of daily prices in the Chinese stock market using 934 842 data records from 2015-01-06 to 2018-01-29, (b2) CDF of daily prices in China stock market, (b3) CDF of 5-minute prices in China stock market from 2019-04-22 to 2019-04-26, (b4) CDF of 10-minute prices in China stock market.

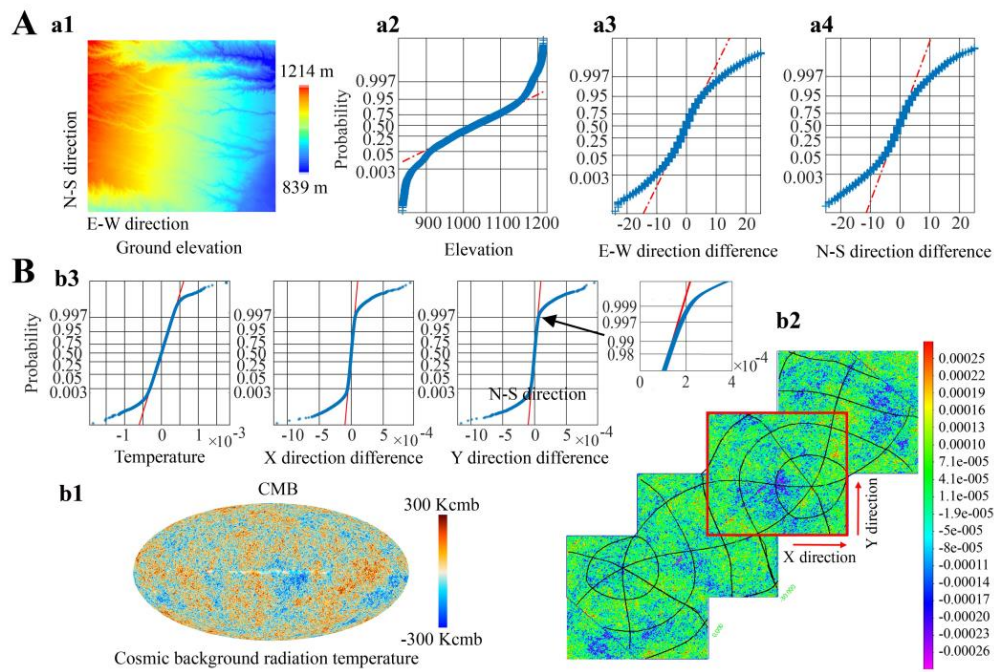

**Fig.9 Spatial distribution of cosmic background radiation and terrestrial elevation. (A)** terrestrial elevation: (a1) ground elevation, (a2) CDF of ground elevation, (a3, a4) CDF of elevation first-order difference in west-east and north-south directions. (**B**) cosmic background radiation: (b1) sky map of CMB radiation temperature (*23*), (b2) data of CMB radiation temperature (SMICA), (b3) CDF of CMB radiation temperature and CMB radiation temperature first-order difference along two directions.

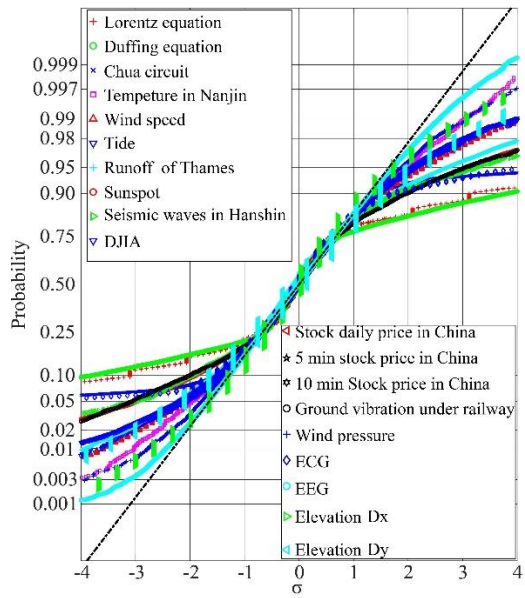

**Fig.10. Normalized cumulative probability functions.**

Table 1. Case information

| Case type | Case name | Sampling period/frequency | Time or space span | Physical background | Sample source |
|---|---|---|---|---|---|
| Classical chaotic dynamical equation | Lorenz system | 0.002 | >100 Cycle | Meteorology | Numerical simulation |
| | Duffing oscillator | 0.02 | >100 Cycle | Mechanical vibration | Numerical simulation |
| | Chua circuit | 0.002 | >100 Cycle | Electro circuit | Numerical simulation |
| Natural phenomena on Earth | Ambient wind speed | 0.1s | 1 hour | Ambient wind speed | Measurement in Yueqing by authors |
| | Tide current | 10min | 30 days | Tide | Measurement in East China Sea |
| | Land runoff | 1 day | 135 years | Hydrology | Web site (UK_NRFA) |
| | Seismic signals | 0.005s~0.01s | 30~60 s | Earthquake | Collection by authors |
| | Sunspot | 1 day | 36 years | Astronomy | Web site (NASA) |
| Social behaviors | Dow Jones stock index (DJIA) | 1 day | 113 years | Economy | Web site (Collected by authors) |
| | Stock price market (China) | 1 day | 2 years | Economy | Web site (sina) |
| | | 5min | 5 days | Economy | Web site (sina) |
| | | 10min | 10 day | Economy | Web site (sina) |
| Biological signals | Human electrocardiogram (adfecgd) | 1000 Hz | 5 min | Biology | Web site (PhysioBank ATM) |
| | Human electrocardiogram (apean) | 1000 Hz | 5 min | Biology | Web site (PhysioBank ATM) |
| | Dog electrocardiogram | 1000Hz | 195 s | Biology | Laboratory measurements, China |
| | PCG (heart tone) | 333 Hz | 20 min | Biology | Web site (PhysioBank ATM) |
| | Brain waves (patient) | 2048 Hz | 10 min | Biology | Web site collection |
| | Brain waves (normal person) | 2048 Hz | 10 min | Biology | Web site collection |
| | Electroencephalogram (ltrsvp) | 2048 Hz | 10min | Biology | Web site (PhysioBank ATM) |
| Artificial structures and experimental device signals | Ground vibration of high-speed train | 1280Hz | 5min | Mechanical vibration | Field measurement |
| | Wind tunnel experiments | 0.01s | 10 min | Aerodynamics | Laboratory measurements, Japan |
| Spatial differential | Cosmic background radiation | 5 arc min | Entire universe | Cosmology | Web site(ESA) |
| | Earth surface elevation | 150 m | 108 km×108km | Geography | Web site (USAG) |

Table 2. Sample *t*-distribution fittings

| | Case | $\frac{\mu}{\sigma}$ | $\nu$ | Correlation coefficient |
|---|---|---|---|---|
| Chaotic equations | Lorentz system | −1.93E−03 | 2.491 | 0.9819 |
| | Duffing equation | −7.77E−06 | 2.658 | 0.9788 |
| | Chua circuit | 2.10E−06 | 1.114 | 0.9413 |
| Natural phenomenon | Temperature | −8.73E−05 | 5.925 | 0.9994 |
| | Wind speed | 1.96E−04 | 3.760 | 0.9980 |
| | Tide current | 3.11E−04 | 3.317 | 0.9988 |
| | Surface runoff | −2.44E−04 | 3.154 | 0.9905 |
| | Sunspot area | −4.41E−04 | 2.561 | 0.9908 |
| | Earth signals | 1.21E−04 | 2.716 | 0.9936 |
| Social time series | Dow Jones index | 1.61E−02 | 2.995 | 0.9989 |
| | Daily stock price in China | 0.0118 | 2.872 | 0.9600 |
| | 5-minute stock price in China | −0.0565 | 2.805 | 0.996 |
| | 10-minute stock price in China | −0.0475 | 2.741 | 0.9840 |
| Biological signal | Dog electrocardiogram | 1.80E−05 | 2.594 | 0.9386 |
| | Brain waves | −5.92E−04 | 13.153 | 0.9998 |
| Artificial structures | Railway vibration | 1.51E−05 | 3.350 | 0.9923 |
| | Wind tunnel pressure | −8.52E−04 | 6.220 | 0.9996 |
| Earth surface elevation | East-west direction | 0.00341 | 4.769 | 0.9931 |
| | North-south direction | −0.01598 | 4.732 | 0.9938 |